\documentclass{article}                                     

\input{amssym.def}                                                           
\input{amssym.tex}                                                           
                                                                             
\begin{document}                                                             
\title{On a desingularization of the moduli space of noncommutative tori}

\author{Igor  ~Nikolaev
\footnote{Partially supported 
by NSERC.}}


 \maketitle


\newtheorem{thm}{Theorem}
\newtheorem{lem}{Lemma}
\newtheorem{dfn}{Definition}
\newtheorem{rmk}{Remark}
\newtheorem{cor}{Corollary}
\newtheorem{prp}{Proposition}
\newtheorem{exm}{Example}
\newtheorem{prb}{Problem}

\newcommand{\N}{{\Bbb N}}
\newcommand{\F}{{\cal F}}
\newcommand{\R}{{\Bbb R}}
\newcommand{\Z}{{\Bbb Z}}
\newcommand{\C}{{\Bbb C}}

\begin{abstract}
It is shown that the moduli space of the
noncommutative tori ${\Bbb A}_{\theta}$ admits  a natural
desingularization by the group $Ext~ ({\Bbb A}_{\theta},{\Bbb A}_{\theta})$.
Namely, we prove  that  the moduli space of pairs $({\Bbb A}_{\theta},  Ext~ ({\Bbb A}_{\theta},{\Bbb A}_{\theta}))$
is homeomorphic to a punctured  two-dimensional sphere.  The proof is based on  
a correspondence (a covariant functor) between the complex and noncommutative tori.

\vspace{7mm}


\vspace{5mm}
{\it AMS (MOS) Subj. Class.:  14H52, 46L85}
\end{abstract}

\section{Introduction}
{\bf A.} Let $0<\theta< 1$ be an irrational number,  whose
continued fraction has the form  $\theta=[a_0, a_1, a_2,\dots]$.
Consider an $AF$-algebra ${\Bbb A}_{\theta}$ given by the following Bratteli diagram:

\begin{figure}[here]
\begin{picture}(300,60)(0,0)
\put(110,30){\circle{3}}
\put(120,20){\circle{3}}
\put(140,20){\circle{3}}
\put(160,20){\circle{3}}
\put(120,40){\circle{3}}
\put(140,40){\circle{3}}
\put(160,40){\circle{3}}

\put(110,30){\line(1,1){10}}
\put(110,30){\line(1,-1){10}}
\put(120,42){\line(1,0){20}}
\put(120,40){\line(1,0){20}}
\put(120,38){\line(1,0){20}}
\put(120,40){\line(1,-1){20}}
\put(120,20){\line(1,1){20}}
\put(140,41){\line(1,0){20}}
\put(140,39){\line(1,0){20}}
\put(140,40){\line(1,-1){20}}
\put(140,20){\line(1,1){20}}

\put(180,20){$\dots$}
\put(180,40){$\dots$}

\put(125,52){$a_0$}
\put(145,52){$a_1$}

\end{picture}

\caption{The Bratteli diagram of the  $AF$-algebra  ${\Bbb A}_{\theta}$.}
\end{figure}
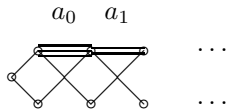

\noindent
where $a_i$ indicate the number of the edges in the upper row of the diagram.  
With a moderate abuse of the terminology,  we shall call  ${\Bbb A}_{\theta}$ 
a  {\it noncommutative torus}.  
(Note that a standard definition of the  noncommutative torus -- a universal $C^*$-algebra
generated by the  unitaries $u,v$ satisfying the commutation relation
$vu=e^{2\pi i}uv$ -- is not an $AF$-algebra. However, the two objects
are isomorphic at the level of their dimension groups \cite{PiVo1}, \cite{Rie1}.)

\medskip
{\bf B.} Recall that the noncommutative tori ${\Bbb A}_{\theta}, {\Bbb A}_{\theta'}$ are said
to be {\it stably isomorphic},  whenever ${\Bbb A}_{\theta}\otimes {\cal K}
\cong {\Bbb A}_{\theta'}\otimes {\cal  K}$,
where ${\cal K}$ is the $AF$-algebra of the compact operators. 
It is well known that the $AF$-algebras  ${\Bbb A}_{\theta}, {\Bbb A}_{\theta'}$ are stably
isomorphic, if and only if,  $\theta'\equiv \theta~mod~SL(2, {\Bbb Z})$,
i.e. $\theta'= (a\theta +b)~/~ (c\theta +d)$,  where 
$a,b,c,d\in {\Bbb Z}$ and $ad-bc=1$ \cite{EfSh1}. 
It is easy to see that the stable isomorphy is an equivalence relation, 
which splits the set $\{{\Bbb A}_{\theta}~|~0<\theta< 1, ~\theta\in {\Bbb R}-{\Bbb Q}\}$ 
into the disjoint equivalence classes. By ${\cal M}$ we shall understand a collection
of such classes, or the ``moduli space'' of  the noncommutative  tori.  An examination
of ${\cal M}$ as a topological space (with the topology
induced by ${\Bbb R}$) shows that the points of ${\cal M}$
have no disjoint neighborhoods, since each orbit  
$\{\theta'\in {\Bbb R}~|~\theta'\equiv \theta~mod~SL(2, {\Bbb Z})\}$
is dense in the real line ${\Bbb R}$. A question arises as of how to ``desingularize''
the (non-Hausdorff) moduli space ${\cal M}$.

\medskip
{\bf C.} Let $A, B$ be a pair of the $C^*$-algebras. Recall that an extension of $A$ by $B$
is a $C^*$-algebra $E$ filling  the  short exact sequence $0\to B\to E\to A\to 0$ of the $C^*$-algebras.
If $A$ is a separable nuclear $C^*$-algebra, the $Ext~(A,B)$ is  an additive abelian  group, 
whose group operation is a sum of the two extensions. 
The $Ext~(A,B)$ is a homotopy invariant in the both variables. The extensions $E_1,E_2$ 
are said to be {\it stably equivalent} if there exists  an isomorphism 
$\psi:E_1\otimes {\cal K}\cong E_2\otimes {\cal K}$,  such  that   
$\psi\circ\alpha_1(B\otimes {\cal K})=\alpha_2(B\otimes {\cal K})$,
where $\alpha_i:B\to E_i, i=1,2$ \cite{B}.  We shall further  
restrict to the case $A=B={\Bbb A}_{\theta}$ and study the stable
equivalence classes of the group $Ext~({\Bbb A}_{\theta}, {\Bbb A}_{\theta})$.  
Using the classification results of D.~Handelman \cite{Han1}, it will
develop that the group $Ext~({\Bbb A}_{\theta}, {\Bbb A}_{\theta})\cong 
Hom~(K_0({\Bbb A}_{\theta}), {\Bbb R})\cong {\Bbb R}$. 
Moreover, the $Ext~({\Bbb A}_{\theta}, {\Bbb A}_{\theta})/$ {\sf stable equivalence}
$\cong {\Bbb R}/{\Bbb Z}$.

\medskip
{\bf D.} An objective of the note is to show that the moduli
of the pairs  
\linebreak
$({\Bbb A}_{\theta}, Ext~({\Bbb A}_{\theta}, {\Bbb A}_{\theta}))$
under the stable equivalence is no longer a non-Hausdorff topological space, but 
a two-dimensional orbifold  (a punctured sphere).  
To prove this result we shall use the Teichm\"uller space  of a torus
(a space of the complex structures on the torus) 
\cite{HuMa1}.  Namely, Hubbard and Masur established a homeomorphism between the Teichm\"uller
space $T_g$ of a surface of genus $g\ge 1$ and the space of quadratic differentials on it. 
We shall use the homeomorphism to extend the action of the modular group $SL(2, {\Bbb Z})$ from 
the upper half-plane ${\Bbb H}=\{x+iy\in {\Bbb C}~|~y>0\}\cong T_1$ to the
space $({\Bbb A}_{\theta}, Ext~({\Bbb A}_{\theta}, {\Bbb A}_{\theta}))$.
Denote by $\widetilde {\cal M}$ the set of pairs 
$({\Bbb A}_{\theta}, Ext~({\Bbb A}_{\theta}, {\Bbb A}_{\theta}))$
modulo the stable equivalence. One obtains the following (natural) desingularization 
of the moduli space of the noncommutative tori. 
\begin{thm}\label{thm1}
$\widetilde {\cal M}$ is a punctured two-dimensional sphere.  
\end{thm}

\bigskip\noindent
{\sf Acknowledgments.} 
I am grateful to Michael Lamoureux for a valuable help in the preparation
of the manuscript. The referee suggestions are kindly acknowledged
and incorporated in the text.

\section{Proof}
We shall split the proof in two lemmas. The background material
is mostly standard, and we shall recall in passing some important notation and 
ideas. 
\begin{lem}\label{lm1}
$\widetilde {\cal M}$ is a two-dimensional  orbifold.
\end{lem}
{\it Proof of lemma \ref{lm1}.}  
We shall use a standard dictionary existing between the $AF$-algebras and their
dimension groups \cite{E}. Instead of dealing with the $AF$-algebra ${\Bbb A}_{\theta}$,
we shall work with its dimension group 
$G_{\theta}=(G, G^+)$, where $G\cong {\Bbb Z}^2$ is the lattice and 
$G^+=\{(x,y)\in {\Bbb Z}^2~|~x+\theta y\ge 0\}$ is a positive cone of
the lattice. The $G_{\theta}$ is the additive abelian group with an order, which  
defines the $AF$-algebra ${\Bbb A}_{\theta}$ up to a stable isomorphism.

Under the dictionary, the extension problem for the $AF$-algebra ${\Bbb A}_{\theta}$
translates as an extension problem for the dimension groups
$G_{\theta}\to E\to G_{\theta}$ (we omit the zeros in the exact sequence). 
An important result of Handelman establishes  the  intrinsic classification
of the extensions of the simple dimension group by a simple dimension group,
see Theorem III.5 of \cite{Han1}. Let us recall  the classification as
it is exposed in \cite{G}, Theorem 17.5 and Corollary 17.7. We shall adopt
the same notation as in the cited work.

Let $H$ be a dense subgroup of the real line ${\Bbb R}$ and $K$ a non-zero dimension
group. Let $E$ be the abelian group $H\oplus K$, and let $\tau: H\to E$
and $\pi: E\to K$ be a natural injection and projection maps. 
Assume that $f: K\to {\Bbb R}$ is a homomorphism of the dimension groups
\footnote{That is  $f$ preserves the positive cone  of $K$ and ${\Bbb R}$:
~$f(K^+)>0$.}
. Then:  (i) $E$ is a dimension group with the positive cone
\begin{equation}\label{eq1}
E^+_f=\{ (0,0) \} ~\bigcup ~\{(x,y)\in E~|~ y\ge 0 ~\hbox{and} ~x+f(y)>0\},
\end{equation}
which gives an extension $H\buildrel{\tau}\over{\to} (E, E^+_f)
\buildrel{\pi}\over{\to} K$ of $H$ by $K$;
(ii) if $f,f': K\to {\Bbb R}$ are the group homomorphisms, then the extensions
$E_f,E_{f'}$ are equivalent if and only if $(f-f')(K)\subseteq H$.

We have to specialize  the above theorem to the case $H=K=G_{\theta}$. 
It is immediate from (i) that $E\cong {\Bbb Z}^4$. Note that
the group homomorphisms $f: G_{\theta}\to {\Bbb R}$ are bijective with the reals
${\Bbb R}$. Indeed, we have to find all the linear maps $f: {\Bbb R}^2\to {\Bbb R}$,
such that $Ker~f= x+\theta y$. (The last equation follows from the condition  $f(G^+)>0$.)
Such maps have the form $f_t(p)=(p,t), ~p,t\in {\Bbb R}^2$, where $(p,t)$ is
the dot product of the two vectors. Let $t=(t_1,t_2)$. Then 
$f_t(-\theta y,y)= t_1(-\theta y)+t_2y= y(t_2-t_1\theta)=0$ for all $y\in {\Bbb R}$.
Therefore, $t_2=\theta t_1$ and $f_t(x,y)=t_1x+\theta t_1y=t_1(x+\theta y), ~t_1\in {\Bbb R}$. 
Thus,  all linear maps $f:{\Bbb R}^2\to {\Bbb R}$ with $Ker~f=x+\theta y$
are bijective with the reals $t_1\in {\Bbb R}$. 
In other words, $Ext~({\Bbb A}_{\theta}, {\Bbb A}_{\theta})$ and ${\Bbb R}$   
are isomorphic as the additive abelian groups.

Let us find when the two extensions $E,E'$ are equivalent. Since $H=G_{\theta}$
is a subgroup of ${\Bbb R}$, one can write $H={\Bbb Z}+\theta {\Bbb Z}$. 
Let $t,t'$ be the real numbers corresponding to the homomorphisms $f,f'$. 
Then $f(G_{\theta})=t({\Bbb Z}+\theta {\Bbb Z})$ and $f'(G_{\theta})=t'({\Bbb Z}+\theta {\Bbb Z})$.
The condition $(f-f')(K)\subseteq H$ of the item (ii) will take the form
$(t-t')({\Bbb Z}+\theta {\Bbb Z})\subseteq {\Bbb Z}+\theta {\Bbb Z}$. 
One  gets immediately that $t=t'+n, n\in {\Bbb Z}$ as a necessary and sufficient condition
for the last inclusion. In other words, the extensions $E,E'$ are equivalent if and only 
if $t'= t~mod ~{\Bbb Z}$. Thus,  the equivalence classes of $Ext~({\Bbb A}_{\theta}, {\Bbb A}_{\theta})$
are bijective with the factor space ${\Bbb R}/{\Bbb Z}$ (a unit interval).

To finish the proof of lemma \ref{lm1}, let us extend the domain of definition
of $\theta$ from the interval $(0,1)$ to the real line ${\Bbb R}$ by allowing $a_0$ to take on any integer value.
In this way,  one can identify the pairs $({\Bbb A}_{\theta}, Ext~({\Bbb A}_{\theta}, {\Bbb A}_{\theta}))$
with the points of ${\Bbb R}^2$ equipped with the usual euclidean topology. We have seen that
the points $(\theta,t)\sim (\theta',t')\in {\Bbb R}^2$ are equivalent if and only if
$\theta'\equiv\theta ~mod ~SL(2,{\Bbb Z})$ and $t'\equiv t ~mod~{\Bbb Z}$. 
Note that the action of the modular group on the second coordinate is always free. Therefore, 
the  points $x,y$ of  the space $\widetilde {\cal M}\cong {\Bbb R}^2/\sim $
admit the disjoint neighborhoods defined, e.g.,  by the open ball of the radius $1/3$ centered  in $x$ and $y$,
respectively. The balls are locally homeomorphic to the euclidean plane, and therefore $\widetilde {\cal M}$ 
is a two-dimensional orbifold. 
$\square$

\bigskip
The lemma \ref{lm1} gives a (partial) desigularization of the space ${\cal M}$.
Indeed, we have seen  that the group $SL(2, {\Bbb Z})\times {\Bbb Z}$ acts in the plane 
$({\Bbb A}_{\theta}, Ext_t~({\Bbb A}_{\theta}, {\Bbb A}_{\theta}))$
by the formula $(\theta,t)\mapsto ({a\theta+b\over c\theta +d}, t+n)$,  where $ad-bc=1$ and $a,b,c,d,n\in {\Bbb Z}$. 
However, the last formula does not specify  the action on  the 
parameter plane $(\theta,t)$ of the modular group  $SL(2,{\Bbb Z})$ alone, 
since the function  $n=n(a,b,c,d)$ is unknown. To find how the integer $n$ depends on 
the integers $a,b,c,d$,  we would need a special construction which involves a correspondence 
(a covariant functor) between  the complex and noncommutative tori. 
Such a construction will be given in the next paragraph and is encapsulated in  
the following lemma.  
\begin{lem}\label{lm2}
There exists a homeomorphism $h: \widetilde {\cal M}
\to {\Bbb H}~/~SL(2,{\Bbb Z}) $, where  ${\Bbb H}=\{x+iy\in {\Bbb C}~|~y>0\}$ is endowed
with hyperbolic metric.
\end{lem}
{\it Proof of lemma \ref{lm2}.}
Let $X$ be a topological surface of genus $g\ge 0$. The Teichm\"uller space
$T_g$ of $X$ consists of the equivalence classes of the complex structures on $X$. 
The $T_g$ is an open ball  of the (real) dimension $6g-6$ if $g\ge 2$ and $2g$ if $g=0,1$. 
By $Mod~X$ we designate a group of the orientation preserving diffeomorphisms
of $X$ modulo the trivial ones. The points $S,S'\in T_g$ are equivalent if there
exists a conformal map $f\in Mod~X$ such that $S'=f(S)$. The moduli 
of conformal equivalence is denoted by ${\cal M}= T_g/~Mod~X$. 
The ${\cal M}$ is a (classical) moduli space,  which  dates back to
Riemann.

Let $S\in T_g$ be a Riemann surface thought as a point in the Teichm\"uller
space, and let $H^0(S, \Omega^{\otimes 2})$ be the space of the holomorphic 
quadratic forms on $S$. The fundamental theorem of Hubbard and Masur
says that there exists a homeomorphism $h_S: H^0(S, \Omega^{\otimes 2})\to T_g$
\cite{HuMa1}, p.224. The  $H^0(S, \Omega^{\otimes 2})$ is a real vector
space of the dimension $6g-6$, where $g\ge2$. It has been shown in the above cited
work,  that $H^0(S, \Omega^{\otimes 2})\cong Hom~(H_1(\tilde X, \tilde \Gamma)^-; {\Bbb R})$
defined by the formula
\begin{equation}\label{eq2}
\omega\longmapsto \left(\gamma\mapsto Im~\int_{\gamma}\omega\right),  
\end{equation}
where $H_1(\tilde X, \tilde \Gamma)^-$ is the odd part in the homology
of a double cover $\tilde X$ of $X$ ramified at zero points $\tilde\Gamma$ 
of the odd multiplicity of the quadratic form \cite{HuMa1}, p.232. 
(The symbolic and formulas will simplify as we come to the complex torus
-- our principal case.) It has been proved that 
$H_1(\tilde X, \tilde \Gamma)^-\cong {\Bbb Z}^{6g-6}$.

Let $X=T^2$, i.e. $g=1$. In this case each  quadratic differential forms
is the  square of a holomorphic abelian form (a one-form), i.e. $H^0(S, \Omega^{\otimes 2})
=H^0(S, \Omega)$.  Therefore $\tilde X=X=T^2$,  $\tilde\Gamma=\emptyset$
and $H_1(\tilde X, \tilde \Gamma)^-=H_1(T^2)\cong {\Bbb Z}^2$.
In other words, one gets a homeomorphism $h_S: Hom~({\Bbb Z}^2, {\Bbb R})\to T_1$.   
As we have seen earlier, $Hom~({\Bbb Z}^2, {\Bbb R})=\{t_1{\Bbb Z}+t_2{\Bbb Z}~|~t_1,t_2\in {\Bbb R}\}=$
$\{t({\Bbb Z}+\theta{\Bbb Z})~|~\theta, t\in {\Bbb R}\}$, where $t=t_1, \theta=t_2 / t_1$. 
On the other hand, the Teichm\"uller space $T_1\cong {\Bbb H}$, where 
${\Bbb H}=\{\tau =x+iy\in {\Bbb C}~|~ y>0\}$ is a (Lobachevsky) upper half-plane and $\tau$ is 
a modulus of the complex torus ${\Bbb C}/({\Bbb Z}+\tau {\Bbb Z})$ \cite{S}, pp. 6-14.  
Thus,  we have a homeomorphism $h_S: ({\Bbb A}_{\theta}, Ext_t~({\Bbb A}_{\theta}, {\Bbb A}_{\theta}))
\to {\Bbb H}$.

Let us show that $h_S$ is equivariant in the first coordinate
with respect to the action of  $Mod~(T^2)\cong SL(2,{\Bbb Z})$,
i.e. $\tau'\equiv\tau ~mod~SL(2,{\Bbb Z})$ if and only if $\theta'\equiv\theta
~mod~SL(2,{\Bbb Z})$. Indeed, since $Hom~(H_1(T^2); {\Bbb R})\cong  {\Bbb H}$
the modular group $SL(2,{\Bbb Z})$ acts on the right-hand side by the formula
$\tau\mapsto (a\tau +b)/(c\tau +d)$ and on the left-hand side by a linear 
transformation  $p_1\mapsto ap_1+bp_2, ~p_2\mapsto cp_1+ dp_2$,
where $p=(p_1,p_2)\in H_1(T^2)$ and $ad-bc=1$.  The $f_p(t)=p_1t_1+p_2t_2$
will become $f_p(t')= t_1(ap_1+bp_2)+t_2(cp_1+dp_2)= p_1t_1'+p_2t_2'$,
where $t_1'=at_1+ct_2$ and $t_2'=bt_1+dt_2$. Therefore  $\theta=t_2/~t_1$
goes to $\theta'=t_2'/t_1'=(b+d\theta)/(a+c\theta)$ and 
$\theta'\equiv\theta ~mod~SL(2,{\Bbb Z})$. (The `only if' part of the 
statement is obtained likewise by an inversion of the formulas.)

\begin{figure}[here]
\begin{picture}(300,150)(0,0)

\put(140,40){\oval(40,40)[t]}
\put(160,40){\oval(40,40)[t]}
\put(180,40){\oval(40,40)[t]}
\put(200,40){\oval(40,40)[t]}
\put(120,40){\oval(40,40)[t]}

\put(130,57){\line(0,1){60}}
\put(110,57){\line(0,1){60}}
\put(150,57){\line(0,1){60}}
\put(170,57){\line(0,1){60}}
\put(190,57){\line(0,1){60}}
\put(210,57){\line(0,1){60}}

\put(60,40){\line(1,0){200}}

\put(157,90){$F$}

\put(157,25){$0$}
\put(132,25){$-1$}
\put(177,25){$1$}
\put(110,25){$-2$}
\put(197,25){$2$}

\end{picture}

\caption{${\Bbb H}$ and the fundamental region $F$.}
\end{figure}
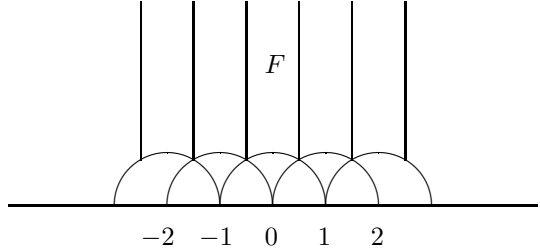

Recall that the Lobachevsky plane ${\Bbb H}=\{x+iy\in {\Bbb C}~|~y>0\}$ 
carries  a  hyperbolic metric $ds=|dz|/y$ such that $SL(2,{\Bbb Z})$
acts on it  by the isometries (linear-fractional transformations). 
The tesselation of ${\Bbb H}$ by the fundamental regions is shown in the Fig.2. 
 Let $\tau'={a\tau +b\over c\tau +d}=T(\tau), ~\tau\in {\Bbb H}$. 
The number $n=n(a,b,c,d)\in {\Bbb Z}$ we shall call a {\it height} of the transformation 
$T\in SL(2,{\Bbb Z})$ if $n$ is equal to the  number of intersections of 
the vertical segment $Im~(\tau'-\tau)$ issued from $\tau$ with the 
lines of tiling ${\Bbb H}/SL(2,{\Bbb Z})$.
(In other words, $n$ shows how many fundamental regions
apart $\tau$ and $\tau'$ are in the vertical direction.) We leave it to the reader
to verify that height $n$ does not depend on a particular choice of $\tau$
in the fundamental region or the fundamental region itself being function of
the transformation $T$ only.

Let us now  define an action of the modular  group 
 on $({\Bbb A}_{\theta}, Ext_t~({\Bbb A}_{\theta}, {\Bbb A}_{\theta}))$.
The action is given by the formula $(\theta,t)\mapsto ({a\theta+b\over c\theta +d}, t+n)$,
where $n=n(a,b,c,d)$ is the  height of the transformation 
$T=T(a,b,c,d)$. Under the homeomorphism $h_S$,  the tesselation of ${\Bbb H}$ 
maps  into a tesselation of the plane $(\theta,t)$.  As we have shown earlier, the  action of
the modular group $SL(2,{\Bbb Z})$ on ${\Bbb H}$ is equivariant with the action 
on $(\theta,t)$.  On the other hand, it is known that ${\Bbb H}/SL(2,{\Bbb Z})$
is a punctured two-dimensional sphere \cite{S}, p.15. 
The lemma \ref{lm2} and theorem \ref{thm1} follow.
$\square$.

\section{Remarks}
Let $E_{\tau}={\Bbb C}/({\Bbb Z}+\tau {\Bbb Z})$ be a complex torus
and $h_S(E_{\tau})=({\Bbb A}_{\theta}, Ext_t~({\Bbb A}_{\theta}, {\Bbb A}_{\theta}))$
its image under the homeomorphism $h_S$. 
Let us call the respective reals $\theta=\theta(\tau)$ and $t=t(\tau)$ a {\it projective curvature}
and an {\it area} of the complex torus $E_{\tau}$.
The projective curvature of the complex tori with a non-trivial group of endomorphisms (complex
multiplication) is a quadratic irrationality. In the latter case, the noncommutative torus
is said to have a {\it real multiplication}.  The noncommutative tori with real multiplication
can be used to  construct the abelian extensions of the real quadratic number fields, 
as  it was suggested  by Yu.~Manin \cite{Man1}.   It seems interesting and challenging  at this point 
to write  a  formula for the projective curvature
and the area as a function of the complex modulus $\tau$.  
It is likely that the  functions  will be of the class $C^0$. 
\begin{prb}
Find a formula (if any) for the functions $\theta(\tau)$ and $t(\tau)$. 
\end{prb}



\vskip1cm

\textsc{The Fields Institute for Mathematical Sciences, Toronto, ON, Canada,  
E-mail:} {\sf igor.v.nikolaev@gmail.com}

\end{document}